\documentclass[12pt,reqno]{amsart}

\usepackage{amssymb, amsmath, amsfonts, latexsym}
\usepackage{enumerate}
\usepackage{color}
\newtheorem{thm}{Theorem}[]
\newtheorem{lem}{Lemma}[section]

\newtheorem{rmk}{Remark}[section]
\theoremstyle{definition}

\numberwithin{equation}{section} \theoremstyle{remark}

\title[Maximum Statistic of standard normal]{\bf Revisit on the convergence rate of normal extremes.}

\author{Liangzhen Lei}
\address{Liangzhen LEI\\ School of Mathematical Science, Capital Normal University, 100048, Beijing, China.} 
\email{leiliangzhen@cnu.edu.cn}

\author{Y{utao} Ma}
\address{Yutao MA\\ School of Mathematical Sciences $\&$ Laboratory  of Mathematics $\&$ Complex Systems of Ministry of Education, Beijing Normal University, 100875 Beijing, China.} 
\thanks{The research of Yutao Ma was supported in part by NSFC 12171038 and 985 Projects.}
\email{mayt@bnu.edu.cn}

\author{Bingjie Tian}
\address{Bingjie Tian\\ School of Mathematical Sciences $\&$ Laboratory of Mathematics and Complex Systems of Ministry of Education, Beijing Normal University, 100875 Beijing, China.}
\email{201711030543@mail.bnu.edu.cn}

\begin{document}
	\maketitle
	
\begin{abstract}
		Let \((X_i)_{1 \le i \le n}\) be independent and identically distributed (i.i.d.) standard Gaussian random variables, and let \(X_{(n)} = \max_{1 \le i \le n} X_i\) denote the maximum order statistic. In extreme value theory, it is well known that the linearly normalized maximum  
\[
Y_n = a_n\bigl(X_{(n)} - b_n\bigr)
\]  
converges weakly to the standard Gumbel distribution \(\Lambda\) as \(n \to \infty\), for suitable scaling and centering constants \(a_n>0\) and \(b_n\).  

In this note, we work with the specific choice  
\[
a_n = \sqrt{2\log n}, \qquad 
b_n = \sqrt{2\log n} - \frac{\log\log n + \log(4\pi)}{2\sqrt{2\log n}},
\]  
and establish the exact rates of convergence of \(Y_n\) to \(\Lambda\) under several probability distances: the Berry-Esseen bound, the \(W_1\) (Wasserstein-1) distance, the total variation distance, the Kullback-Leibler divergence, and the Fisher information. We also illustrate how the orders of convergence depend on the selection of the sequences \(a_n\) and \(b_n\).
   
	\end{abstract} 
	
	
	{\bf AMS Classification Subjects 2020:} 60G70, 60F05, 62G32.

	\section{Introduction} 
	Let $(X_i)_{1 \le i \le n}$ be independent and identically distributed (i.i.d.) standard Gaussian random variables, and denote by $X_{(n)} = \max_{1 \le i \le n} X_i$ the maximum order statistic. It is well-known in extreme value theory that the linearly normalized maximum  
$
Y_n = a_n(X_{(n)} - b_n),
$  
converges weakly to the standard Gumbel distribution $\Lambda$ as $n \to \infty$, where $a_n > 0$ and $b_n$ are appropriate scaling and centering constants. The cumulative distribution function of the limiting distribution is given by $\Lambda(x) = e^{-e^{-x}}$ for $x \in \mathbb{R}.$ 
 
The appropriate choice of the normalizing sequences $a_n$ and $b_n$ has a rich historical background, with foundational contributions from several key works in extreme value theory. The problem was first studied by Fisher and Tippett (1928) in their pioneering work on limiting forms of extreme order statistics. Later, Gnedenko (1943) rigorously established necessary and sufficient conditions for the convergence of maxima to each of the three extreme value distributions, with the Gaussian case falling into the Gumbel domain of attraction. Cram\'er (1946) provided further refinements, deriving explicit expressions for the normalizing constants in the Gaussian case.  

The standard choices for the normalization sequences are:  
\begin{equation}\label{classi}
a_n = \sqrt{2 \log n} \quad \text{and} \quad b_n = \sqrt{2 \log n} - \frac{\log \sqrt{4\pi\log n}}{\sqrt{2 \log n}},
\end{equation}   
which ensure the weak convergence of \(Y_n\) to the standard Gumbel distribution. These constants emerge from careful analysis of the tail behavior of the Gaussian distribution and the theory of extreme value convergence.  
	
Every convergence result is accompanied by the question of rate of convergence. 
The rate at which the normalized $X_{(n)}$ converges to the Gumbel distribution was first rigorously studied by Hall (\cite{Hall}, 1979). He established that when the normalization sequences are chosen as \( b_n \) being the solution to \( 2\pi b_n^2 \exp(b_n^2) = n^2 \),   
 there exist absolute constants \( c_1, c_2 > 0 \) (independent of \( n \)) such that:  

\begin{equation}\label{hallberry}
\frac{c_1}{\log n} \leq \sup_{x \in \mathbb{R}} \left| \mathbb{P}\left( X_{(n)} \leq b_n + \frac{x}{b_n} \right) - e^{-e^{-x}} \right| \leq \frac{c_2}{\log n},  
\end{equation}  
with \( c_2 \leq 3 \) for sufficiently large \( n \). Hall further proved that the \( \frac{1}{\log n} \) convergence rate is {\it sharp}-no faster uniform convergence can be achieved by any alternative choice of norming sequences. Interestingly,
Leadbetter, Lindgrenand Rootz\'en \cite[(2.4.8.)]{Leadbetter1983} proved that when $a_n$ and $b_n$ are chosen as in \eqref{classi}, 
\begin{equation}\label{pointwise}
\mathbb{P}(X_{(n)}\le b_n+\frac{x}{a_n})=e^{-e^{-x}}+e^{-e^{-x}-x}\frac{(\log\log n)^2}{16 \log n}(1+o(1))	
\end{equation}
for any $x\in \mathbb{R}$ fixed.  This highlights the crucial role of norming constant selection in optimizing the speed of convergence. 

Our primary motivations stem from two key considerations. First, we aim to refine the Berry-Esseen bound given in \eqref{hallberry} by Hall \cite{Hall}, with the specific goal of deriving a precise asymptotic expression for the convergence rate. Second, a series of studies conducted by the second author in collaboration with others have explored the exact order of convergence for extreme eigenvalues of non-Hermitian random matrices to Gumbel distributions, 
under both the Berry-Esseen bound and the $W_1$ Wasserstein distance (\cite{MaMeng2025, MaWang2025}). These investigations have inspired us  to revisit this classical framework.

As a natural parallel inquiry, we are particularly focused on determining the exact convergence rates under various other probability metrics. Commonly used metrics in this context include the $W_1$ Wasserstein distance, total variation distance, Kullback-Leibler divergence, and Fisher information distance

To facilitate our discussion, we first recall the formal definitions of these metrics. Let $\mu$ and $\nu$ be two probability distributions on $\mathbb{R}$ with density functions $f$ and $g$, respectively. The distances and divergences are defined as follows (see \cite{Villani2000}).  	
	\begin{itemize}
	\item The $W_1$ Wasserstein distance of $\mu$ and $\nu$ are given by 
	$$W_1(\mu, \nu)=\int_{\mathbb{R}}|F_{\mu}(x)-F_{\nu}(x)|dx.$$ 
	\item The total variation distance is $$\|\mu-\nu\|_{\rm TV}=\int_{\mathbb{R}}|f(x)-g(x)|dx.$$
	\item The Kullback-Leibler divergence is given by 
	$$D_{\rm KL}(\mu|\nu)=\int_{\mathbb{R}} f(x)\log\frac{f(x)}{g(x)} dx.$$
	\item The Fisher information is 
	$$I(\mu|\nu)=\int_{\mathbb{R}}|(\sqrt{f/g})'(x)|^2 g(x)dx,$$	
	which could be simplified as 
	\begin{equation}\label{Fishervar}
		I(\mu|\nu)=\frac{1}{4}\int_{\mathbb{R}} \frac{g^2(x)}{f^2(x)}[(\frac{f}{g})'(x)]^2 f(x)dx=\frac{1}{4}\int_{\mathbb{R}}[(\log \frac{f}{g})'(x)]^2 f(x) dx.
	\end{equation}
	\end{itemize}

	We now establish the exact asymptotic rates for the Berry-Esseen bound, along with the four aforementioned distance metrics between the distribution of $Y_n$ and the Gumbel distribution $\Lambda$.

	\begin{thm}\label{main} Let $\{X_i\}_{i\ge 1}$ be i.i.d. standard normal and $X_{(n)}=\max_{1\le i\le n} X_i.$ Let $$a_n = \sqrt{2 \log n} \quad \text{and} \quad b_n = \sqrt{2 \log n} - \frac{\log (\sqrt{4\pi\log n}) }{\sqrt{2\log n}}$$ and set $Y_n=a_n(X_{(n)}-b_n).$ Let $\Lambda$ be the Gumbel distribution and let $\mathcal{L}(Y_n)$  be the distribution of $Y_n.$  We have the following assertions as $n$ large enough.   
		$$\aligned \sup_{x\in \mathbb{R}}|\mathbb{P}(Y_n\le x)-\Lambda(x)|&=\frac{(\log\log n)^2}{16e\log n}(1+o(1));\\
			W_1(\mathcal{L}(Y_n), \Lambda)&=\frac{(\log\log n)^2}{16\log n}(1+o(1)); \\
			\|\mathcal{L}(Y_n)-\Lambda\|_{\rm TV}&=\frac{(\log\log n)^2}{8e\log n}(1+o(1));\\
			D_{\rm KL}(\mathcal{L}(Y_n)|\Lambda)&=\frac{(\log\log n)^4}{512(\log n)^2}(1+o(1)); \\
			I(\mathcal{L}(Y_n)|\Lambda)&=\frac{(\log\log n)^4}{512(\log n)^2}(1+o(1)).
			 \endaligned $$
	\end{thm}

	In the second section, we establish key lemmas, which include precise asymptotic results for both the distribution function and the density function of $Y_n$. The third section is dedicated to prove Theorem \ref{main}. Finally, we discuss how the choices of $a_n$ and $b_n$ 
  affect the convergence rate in the forth section.
	
	Throughout this paper, we use the following notation:
	\begin{itemize}
		\item $\ell_1(n) =  \frac14\log\log n, \quad \ell_2(n) =c_n-1.$
		\item $\phi(x)$ and $\Phi(x)$ denote the density and cumulative distribution function of  the standard normal, respectively and $\Psi(x) = 1 - \Phi(x)$ represents the complementary cumulative distribution function of the standard normal.
		\item For two positive sequence $z_n$ and $t_n,$ the notation $z_n\lesssim t_n$(equivalently $t_n\gtrsim z_n$) means that  there exists a constant $c>0$ such that $z_n\le c t_n$ and $z_n\asymp t_n$ if $z_n\lesssim t_n$ and $t_n\lesssim z_n.$ 
		\item $z_n=O(t_n)$ means $\varlimsup\limits_{n\to\infty}\frac{|z_n|}{t_n}\le c$ for some $c>0.$
		
	\end{itemize} 
	
	\section{Premililaries} 
In this section, we always use the following notations. 
$$a_n:= \sqrt{2 \log n},  \quad c_n = \log(\sqrt{2\pi} a_n) \quad \text{and} \quad t_n(x)=x-c_n.$$	 
Now $$\ell_1(n) =  \frac14\log\log n, \quad \ell_2(n) =(\log n)^{1/4}$$ and then  
\begin{equation}\label{tnxregime}
	|t_n(x)|\lesssim (\log n)^{1/4}
\end{equation}
uniformly on $x\in [-\ell_1(n), \ell_2(n)].$ The inequalities  \eqref{tnxregime} will be used frequently.

The metrics that we consider are fundamentally connected to either the distribution function or the density function of $Y_n$. To establish our results, we first derive precise asymptotic expressions for these functions. Let $F_n$ denote the distribution function and $f_n$ the density function of $Y_n,$ respectively. By definition, it holds that 
\begin{equation}\label{distri} F_n(x)=\mathbb{P}(Y_n\le x)=\mathbb{P}(X_{(n)}\le a_n+\frac{x-c_n}{a_n})=\Phi^n(a_n+\frac{t_n(x)}{a_n})\end{equation} 
and then 
\begin{equation}\label{densi} f_n(x)=n a_n^{-1}\varphi\left(a_n + \frac{t_n(x)}{a_n}\right)\Phi^{n-1}\left(a_n + \frac{t_n(x)}{a_n}\right). \end{equation} 

To lay the groundwork for our main theorem, we establish several key technical lemmas. These results provide precise asymptotic approximations for the distribution function \( F_n(x) \) and the density function \( f_n(x) \) for $x\in [-\ell_1(n), \ell_2(n)]$, which play a central role in our subsequent analysis.

	\begin{lem}\label{psilem} Let $a_n$ and $t_n(x)$ be given as above. We have
		\begin{equation*}
			\Psi\big(a_n+\frac{t_n(x)}{a_n}\big)=\frac{1}{n}e^{-x}\big(1-\frac{t_n^2(x)+2t_n(x)+2}{4\log n}\big)(1+O(\frac{t_n^4(x)+1}{(\log n)^{2}}))
		\end{equation*}
		uniformly on $x\in [-\ell_1(n), \ell_2(n)].$
	\end{lem}
	
	\begin{proof}
		The inequalities \eqref{tnxregime} makes sure  \begin{equation}\label{anexpre}
			a_n+\frac{t_n(x)}{a_n}=\sqrt{2\log n}(1+\frac{t_n(x)}{2\log n})\asymp \sqrt{\log n}\gg 1\end{equation} uniformly for $x \in [-\ell_1(n), \ell_2(n)]$ for sufficiently large  $n.$ The well-known Mills ratio derives 
		\begin{equation}\label{Millsratio}
			\Psi\big(a_n+\frac{t_n(x)}{a_n}\big)=\frac{\exp(-\frac12(a_n+\frac{t_n(x)}{a_n})^2)}{\sqrt{2\pi}\big(a_n+\frac{t_n(x)}{a_n}\big)}(1-(a_n+\frac{t_n(x)}{a_n})^{-2})\big(1+O\big(\frac{1}{a_n^4}\big)\big).
		\end{equation}
		The fact $c_n=\log(\sqrt{2\pi} a_n), a_n=\sqrt{2\log n}$ and \eqref{anexpre} lead   
		\begin{equation}\label{somepoint}\aligned \exp(-\frac{1}{2}(a_n+\frac{t_n(x)}{a_n})^2)&=\exp(-\log n-x+\log(\sqrt{2\pi}a_n)-\frac{t_n^2(x)}{4\log n})\\
		&=\frac{\sqrt{2\pi} a_n}{n}e^{-x} \exp(-\frac{t_n^2(x)}{4\log n}),\endaligned \end{equation}
		whence 
		\begin{equation}\label{Psiexp}
			\Psi\big(a_n+\frac{t_n(x)}{a_n}\big)=\frac{1}{n} e^{-x} \frac{e^{-\frac{t_n^2(x)}{4 \log n}}}{1+\frac{t_n(x)}{2 \log n}}\big(1-\frac{1}{2\log n}(1+\frac{t_n(x)}{2\log n})^{-2}\big)(1+O(\frac{1}{(\log n)^2})).
		\end{equation}
		Applying  corresponding  Taylor expansions while keeping in mind $|t_n(x)|\lesssim (\log n)^{1/4},$ we obtain
			\begin{equation}\label{expansion}\aligned 
				\exp(-\frac{t_n^2(x)}{4 \log n})&=1-\frac{t_n^2(x)}{4 \log n}+O\big(\frac{t_n^4(x)+ 1}{(\log n)^2}\big);\\
			(1+\frac{t_n(x)}{2 \log n})^{-1}&=1-\frac{t_n(x)}{2 \log n}+O\big(\frac{t_n^2(x)+1}{(\log n)^2}\big);\\
			1-\frac{1}{2\log n}(1+\frac{t_n(x)}{2\log n})^{-2}&=1-\frac{1}{2\log n}+O(\frac{|t_n(x)|+1}{(\log n)^2}). \endaligned 
		\end{equation}
		Putting these asymptotics back into \eqref{Psiexp},  we pick up the terms up to the order $\frac{1}{\log n}$ to derive 
		\begin{equation*}
			\Psi\left(a_n+\frac{t_n(x)}{a_n}\right)=\frac{1}{n}e^{-x}(1-\frac{t_n^2(x)+2(t_n(x)+1)}{4\log n})(1+O(\frac{t_n^4(x)+1}{(\log n)^{2}})).
		\end{equation*}
		This completes the proof. 
	\end{proof}

	\begin{lem}\label{keyFesti} Let $F_n$ be defined as above. Then,   
		\begin{align*}
			F_n(x)-e^{-e^{-x}}=e^{-e^{-x}} e^{-x}\frac{t_n^2(x)+2t_n(x)+2}{4 \log n}(1+O((\log n)^{-1}))	\end{align*}
		uniformly on $[-\ell_1(n), \ell_2(n)]$ as $n\to\infty.$ 
	\end{lem}

	\begin{proof}
		By definition, we have 
		\begin{equation}\label{Fninitial}
			F_n(x)=\Phi^n\big(a_n+\frac{t_n(x)}{a_n}\big)=\exp\big(n\log (1-\Psi(a_n+\frac{t_n(x)}{a_n}))\big).
		\end{equation}
Lemma \ref{psilem} indicates $$n\Psi(a_n+\frac{t_n(x)}{a_n})=e^{-x}(1-\frac{t_n^2(x)+2t_n(x)+2}{4\log n})(1+O(\frac{t_n^4(x)+1}{(\log n)^{2}}))$$ as well as 
$$n\Psi^2(a_n+\frac{t_n(x)}{a_n})\lesssim \frac{1}{n} e^{-2 x}\le\frac{\exp(2\ell_1(n))}{n}=\frac{(\log n)^{1/2}}{n}. $$
Thereby,  
applying the Taylor expansion on the righthand side of \eqref{Fninitial} and using Lemma \ref{psilem}, we derive 		\begin{align*}
			F_n(x)&=\exp\big(-n\Psi(a_n+\frac{t_n(x)}{a_n})+O(\frac{(\log n)^{1/2}}{n})\big)\\
			&=\exp\left\{-e^{-x}+e^{-x}\frac{t_n^2(x)+2t_n(x)+2}{4\log n}+O(\frac{e^{-x}(t_n^4(x)+1)}{(\log n)^2}+\frac{(\log n)^{1/2}}{n})\right\}.
		\end{align*}
	The imposed constraint on 
$x$
and \eqref{tnxregime} guarantee that the following key inequality holds:
\begin{equation}\label{01}
e^{-x} \frac{t_n^2(x) + 2t_n(x) + 2}{4 \log n} \lesssim \frac{(\log n)^{1/2} e^{\ell_1(n)}}{\log n}=(\log n)^{-1/4}= o(1),
\end{equation}
which serves as the fundamental rationale behind our specific choice of the range for 
$x$. 
	Thus, applying a Taylor expansion for $e^t=1+t+O(t^2)$ for $|t|\ll 1,$ we continue to have 
	\begin{equation}\label{Fnasym}\aligned 
F_n(x)&=e^{-e^{-x}}\Big(1+
 e^{-x} \frac{t_n^2(x) + 2t_n(x) + 2}{4 \log n}+O(e^{-2x}\frac{(t_n^2(x) + 2t_n(x) + 2)^2}{(\log n)^2})\\
 &\quad \quad \quad \quad \quad\quad \quad \quad\quad \quad \quad\quad\quad  +O(\frac{e^{-x}(t_n^4(x)+1)}{(\log n)^2}+\frac{(\log n)^{1/2}}{n})\Big).
\endaligned \end{equation}
The constraint on $x$ makes sure that the two error terms in \eqref{Fnasym}	have order not greater than   
	$e^{-x}\frac{(\log\log n)^4}{(\log n)^{7/4}}$ based on \eqref{tnxregime}, whence   
	\begin{equation}\label{Fnasym1}\aligned
F_n(x)&=e^{-e^{-x}}(1+
 e^{-x} \frac{t_n^2(x) + 2t_n(x) + 2}{4 \log n}+O(\frac{(e^{-x}\vee e^{-2x})(t_n^4(x)+1)}{(\log n)^{2}}))\\
&=e^{-e^{-x}}+\exp(-e^{-x}-x)\frac{t_n^2(x) + 2t_n(x) + 2}{4 \log n}(1+O(\frac{(e^{-x}\vee 1)(t_n^4(x)+1)}{(\log n)^{2}})).
\endaligned
\end{equation}
Considering the function $h(x)=e^{-x}(x-c_n)^4,$ it is straightforward to see 
$f'>0$ on $(c_n, c_n+4)$ and $f'<0$ otherwise. Therefore, 
$$e^{-x}t_n^4(x)\le \max\{h(-\ell_1(n)), h(c_n+4)\}\lesssim \exp(\ell_1(n))(\log\log n)^4\ll \log n.$$ 
Thus, by the expression \eqref{tnxregime}, we have  
	\begin{equation}\label{Oerr}
	(e^{-x}\vee 1)(t_n^4(x)+1)\lesssim \log n	
	\end{equation}
	and putting \eqref{Oerr} back into \eqref{Fnasym1} finishes 
	the proof. 
			\end{proof}

\begin{rmk}
When $x\in\mathbb{R}$ is fixed, the expression \eqref{Fnasym} and the fact $$t_n(x)=c_n(1+o(1))=\frac{1+o(1)}{2}\log \log n$$ derive that   
$$F_n(x)=e^{-e^{-x}}+\exp(-e^{-x}-x)\frac{(\log\log n)^2}{16\log n}(1+o(1)).$$
This is exactly the expression \eqref{pointwise} stated as (2.4.8) in \cite{Leadbetter1983}.	
\end{rmk}

Now we give a precise asymptotic on $f_n.$ 
\begin{lem}\label{densilem} Let $f_n$ be given in \eqref{densi}. We have 
$$f_n(x)=\exp\big(-x-e^{-x}\big) \big(1+e^{-x}\frac{t_n^2(x)+2t_n(x)+2}{4 \log n}-\frac{t_n^2(x)}{4\log n}\big)(1+O((\log n)^{-1}))$$ for $n$ large enough when $x\in [-\ell_1(n), \ell_2(n)].$ 
\end{lem}

 \begin{proof}
 Inserting the expression \eqref{somepoint} into \eqref{densi}, we obtain 
 \begin{equation}\label{halffn} f_n(x)=\exp\big(-x-\frac{t_n^2(x)}{4\log n}\big)(F_n(x))^{\frac{n-1}{n}}.\end{equation}
 We write $$(F_n(x))^{\frac{n-1}{n}}=F_n(x)\exp(-\frac1n\log F_n(x)).$$ Using the first equality in \eqref{Fnasym1}, we see 
$$ -\frac{1}{n}\log F_n(x)=\frac{e^{-x}}{n}+O(e^{-x}\frac{t_n^2(x) +1}{n\log n}+\frac{(e^{-x}\vee e^{-2x})(t_n^4(x)+1)}{n(\log n)^{2}})),$$ 
which ensures 
$$\exp(-\frac1n\log F_n(x))-1\asymp \frac{e^{-x}}{n}\lesssim  \frac{(\log n)^{1/4}}{n}.$$ 
Then, we have 
$$\aligned f_n(x)&=\exp\big(-x-\frac{t_n^2(x)}{4\log n}\big)F_n(x)(1+O(\frac{(\log n)^{1/4}}{n}))\\
&=\exp\big(-x-\frac{t_n^2(x)}{4\log n}-e^{-x}\big)(1+e^{-x}\frac{t_n(x)^2+2t_n(x)+2}{4 \log n})(1+O((\log n)^{-1})).\endaligned $$ Continuing to write $$\exp(-\frac{t_n^2(x)}{4\log n})=1-\frac{t_n^2(x)}{4\log n}+O(\frac{t_n^4(x)+1}{(\log n)^2})=1-\frac{t_n^2(x)}{4\log n}+O((\log n)^{-1}),$$ 
we finally derive 
$$\aligned f_n(x)=\exp(-x-e^{-x})(1+e^{-x}\frac{t_n(x)^2+2t_n(x)+2}{4 \log n}-\frac{t_n^2(x)}{4\log n})(1+O((\log n)^{-1})).\endaligned $$
The proof is completed. 
 \end{proof}

\begin{rmk} In our subsequent analysis of the Kullback-Leibler divergence, a crucial component involves evaluating the expectation $\mathbb[e^{-Y_n}-\frac{t_n^2(Y_n)}{4\log n}].$ Initial attempts using Lemma \ref{densilem}, where we expanded $f_n$ via Taylor's formula with remainder $O(\frac{(t_n^{4}(x)+1}{(\log n)^{2}}),$ proved inadequate for capturing the leading-order behavior of the divergence. This insufficiency becomes apparent when examining the delicate balance between the various terms contributing to the divergence. Replacing \eqref{Millsratio} by 
$$\Psi(t)=\frac{1}{\sqrt{2\pi} t}(1-t^{-2}+3 t^{-4})(1+O(t^{-6}))$$ for $t\gg 1,$ using \eqref{somepoint} and expanding all three left sides of \eqref{expansion} up to terms with denominator $(\log n)^{-3},$ we get  via straightforward calculation that 
$$\aligned &\quad\Psi(a_n+\frac{t_n(x)}{a_n})\\
=&\frac{\exp(-\frac{1}{2}(a_n+\frac{t_n(x)}{a_n})^2)}{\sqrt{2\pi}\big(a_n+\frac{t_n(x)}{a_n}\big)}(1-(a_n+\frac{t_n(x)}{a_n})^{-2}+3(a_n+\frac{t_n(x)}{a_n})^{-4})\big(1+O\big(\frac{1}{a_n^6}\big)\big)\\
=&\frac{\exp(-x-\frac{t_n^2(x)}{4 \log n})}{n(1+\frac{t_n(x)}{2 \log n})}\big(1-\frac{1}{2\log n}(1+\frac{t_n(x)}{2\log n})^{-2}+\frac{3}{4(\log n)^2}(1+\frac{t_n(x)}{2\log n})^{-4}\big)(1+O(\frac{1}{(\log n)^3})) \\
=&\frac{e^{-x}}{n}(1-\frac{t_n^2(x)+2t_n(x)+2}{4\log n}+\frac{t_n^4(x)+4t_n^3(x)+12t_n^2(x)+24t_n(x)+24}{32(\log n)^2})(1+O(\frac{t_n^6(x)+1}{(\log n)^3})).\endaligned $$ 
 Plugging this asymptotic into the expression of $f_n,$ while always extending the Taylor formula up to $O(\frac{t_n^6(x)}{(\log n)^3}),$ we get the a more precise expression of $f_n$ as 
$$\aligned f_n(x)&=\exp(-x-e^{-x})\bigg(1+e^{-x}\frac{t^2_n(x)+2t_n(x)+2}{4 \log n}-\frac{t_n^2(x)}{4\log n}\\
&\quad +\frac{t_n^4(x)}{32(\log n)^2}-\frac{e^{-x}}{32(\log n)^2}\big(3t_n^4(x)+8t_n^3(x)+16t_n^2(x)+24t_n(x)+24\big)\\
&\quad+\frac{e^{-2x}}{32(\log n)^2}(t_n^2(x)+2t_n(x)+2)^2\bigg)(1+O(\frac{(t_n^6(x)+1)(e^{-x}\vee 1\vee e^{-3x})}{(\log n)^{3}})).\endaligned$$
The same analysis as for \eqref{Oerr} leads  
$$(t_n^6(x)+1)(e^{-x}\vee 1\vee e^{-3x})\lesssim (\log n)^2$$ and then we get 
\begin{equation}\label{densifnew}\aligned f_n(x)&=\exp(-x-e^{-x})\bigg(1+e^{-x}\frac{t^2_n(x)+2t_n(x)+2}{4 \log n}-\frac{t_n^2(x)}{4\log n}\\
&\quad +\frac{t_n^4(x)}{32(\log n)^2}-\frac{e^{-x}}{32(\log n)^2}\big(3t_n^4(x)+8t_n^3(x)+16t_n^2(x)+24t_n(x)+24\big)\\
&\quad+\frac{e^{-2x}}{32(\log n)^2}(t_n^2(x)+2t_n(x)+2)^2\bigg)(1+O((\log n)^{-1})).\endaligned\end{equation} The refined expansion reveals that
the terms up to $\frac{t_n^4(x)}{(\log n)^2}$ contribute meaningfully to the expectation, 
the sixth-order expansion captures previously overlooked cancellations, 
and the resulting approximation achieves the tight error bounds needed for our divergence analysis.
This enhanced approximation scheme provides the mathematical resolution necessary to establish the precise asymptotic behavior of the Kullback-Leibler divergence in our setting. 

\end{rmk}

The derivation of precise asymptotic expansions for the density function \( f_n \) allows us to obtain refined moment estimates for \( Y_n \). 
Furthermore, our analysis also yields an enhanced asymptotic expansion for the expectation
$\mathbb{E}[e^{-Y_n}-\frac{t_n^2(Y_n)}{4\log n}],$
which plays a pivotal role in our subsequent analysis of Kullback-Leibler divergence. Even we don't need the following refined expression for $\mathbb{E}Y_n,$ we still present it to illustrate how to get the precise asymptotic of $\mathbb{E}[e^{-Y_n}-\frac{t_n^2(Y_n)}{4\log n}].$ 

	\begin{lem}\label{expec} The following asymptotic formulas hold
			\begin{align*}
			\mathbb{E}[Y_n]&=\gamma-\frac{c_n^2}{4\log n}+\frac{\gamma c_n}{2\log n}+O((\log n)^{-1});\\
			\mathbb{E}(e^{-Y_n}-\frac{t_n^2(Y_n)}{4\log n})&=1+\frac{c_n^4}{32(\log n)^2}(1+o(1))
		\end{align*}
		for sufficiently large $n.$ Here $\gamma$ denotes the Euler-Mascheroni constant.
	\end{lem}

	\begin{proof}
			We first present the detailed derivation for $\mathbb{E}[Y_n]$ as an illustrative example.
 By definition, 	
 $$\mathbb{E}[Y_n]=\int_{\mathbb{R}} x f_n(x)d x$$	and we cut the integration interval $\mathbb{R}$ into three parts, according to Lemma \ref{densi}, as  	\begin{align*}
				\mathbb{E}[Y_n]=\left(\int_{-\infty}^{-\ell_1(n)}+\int_{-\ell_1(n)}^{\ell_2(n)}+\int_{\ell_2(n)}^{+\infty}\right) x f_n(x) d x\triangleq& \mathrm{I}+\mathrm{I\!I}+\mathrm{I\!I\!I},
			\end{align*}
			where $\ell_1(n)=\frac14\log\log n, \; \ell_2(n)=(\log n)^{1/4}.$
		We first estimate the term $\mathrm{III}.$	
		Indeed, we use \eqref{halffn} and the fact $\Phi\le 1$ to get  
		\begin{equation}\label{estiiii} \aligned & \int^{+\infty}_{\ell_2(n)} x \exp(-x-\frac{t_n^2(x)}{4\log n})\Phi^{n-1}(a_n+\frac{t_n(x)}{a_n}) dx\\
		&\le \int^{+\infty}_{\ell_2(n)} x \exp(-x)dx=(1+\ell_2(n))e^{-\ell_2(n)}\ll\frac{1}{\log n}.
		\endaligned \end{equation}					
	Meanwhile,	for the integral $\mathrm{I}$, we use \eqref{halffn} and the monotonicity of $\Phi$  to get  
			\begin{equation}\label{estiupi}\aligned &\quad \int_{-\infty}^{-2\log n} (-x) \exp(-x-\frac{t_n^2(x)}{4\log n})\Phi^{n-1}(a_n+\frac{t_n(x)}{a_n}) dx \\
			&=\int^{+\infty}_{2\log n} t \exp(t-\frac{(t+c_n)^2}{4\log n})\Phi^{n-1}(a_n-\frac{t+c_n}{a_n}) d t\\
			&\le \Phi^{n-1}(a_n-\frac{2\log n}{a_n}) \int^{\infty}_{2\log n} t \exp(t-\frac{(t+c_n)^2}{4\log n}) dt\\
			&=2^{-n}\int^{\infty}_{2\log n} t \exp(t-\frac{(t+c_n)^2}{4\log n}) dt. \endaligned \end{equation}	
		We simplify the integral at the right hand side of the last line  in \eqref{estiupi} as 
		$$\int^{+\infty}_{2\log n} t \exp(t-\frac{(t+c_n)^2}{4\log n}) dt=\int^{+\infty}_{2\log n} \exp(-w(t)) dt$$
		by setting $$w(t)=-\log t-t+\frac{(t+c_n)^2}{4\log n}.$$ Now, the derivative of $w$ is 
		$$w'(t)=-1-\frac1t+\frac{t+c_n}{2\log n},$$ which is increasing and then 
		$$w'(t)\ge w'(2\log n)=\frac{c_n-1}{2\log n}>0$$ since $t\ge 2\log n.$ 
		Thus, it holds from the fact $w(+\infty)=+\infty$ that  
			\begin{equation}\label{intei1}\aligned
			\int^{+\infty}_{2\log n} e^{-w(t)} dt&\le \frac{1}{w'(2\log n)}\int^{+\infty}_{2\log n} e^{-w(t)} w'(t) dt\\
			&=\frac{\exp(-w(2\log n))}{w'(2\log n)}\lesssim n (\log n)^{3/2}(\log\log n)^{-1}. 
			\endaligned \end{equation}
			Inserting \eqref{intei1} into \eqref{estiupi}, we get 
			$$\int_{-\infty}^{-2\log n}|x| f_n(x)dx\lesssim n (\log n)^{3/2}(\log\log n)^{-1}\, 2^{-n}\ll(\log n)^{-1}.$$
			Now we treat the remainder of $|\mathrm {I}|,$ which is in fact 
			$\int_{-2\log n}^{-\ell_1(n)} |x| f_n(x) dx.$ Simple algebra and \eqref{Fnasym1} bring directly 
			$$\aligned \int_{-2\log n}^{-\ell_1(n)} |x| f_n(x) dx&\lesssim (\log n) F_n(-\ell_1(n))
			\lesssim \log n\, e^{-(\log n)^{1/4}}\ll (\log n)^{-1}.\endaligned$$  
			Thus,  we get 
			$$|\mathrm I|+\mathrm{I\!I\!I}\ll (\log n)^{-1}.$$
		Now, it is the time to work on the most difficult part $\mathrm{I\!I}.$ Lemma \ref{densilem}  		
			enables us to write 
			\begin{equation}\label{IIsum}\aligned\mathrm{I\!I}&=\int_{-\ell_1(n)}^{\ell_2(n)} xf_n(x) dx\\
			&=(1+O((\log n)^{-1}))\int_{-\ell_1(n)}^{\ell_2(n)} xe^{-x-e^{-x}}(1+e^{-x}\frac{t_n(x)^2+2t_n(x)+2}{4 \log n}-\frac{t_n^2(x)}{4\log n})dx.\endaligned \end{equation}
			As we know, on the one hand 
			$$\int_{-\infty}^{+\infty} x e^{-x-e^{-x}} dx=\gamma.$$ 
			On the other hand, the substitution of variable implies 
			$$\aligned -\int_{-\infty}^{-\ell_1(n)} x e^{-x-e^{-x}} dx&=\int^{+\infty}_{(\log n)^{1/4}} e^{-t} \log t dt\\
			&=\frac14(\log\log  n)e^{-(\log n)^{1/4}}+\int^{+\infty}_{(\log n)^{1/4}}\frac1t e^{-t} dt\\
			&\lesssim (\log\log  n)e^{-(\log n)^{1/4}}  \endaligned $$
						and similarly
						$$\aligned\int^{+\infty}_{\ell_2(n)} x e^{-x-e^{-x}} dx&=\int_0^{e^{-(\log n)^{1/4}}}(-\log t) e^{-t} dt\\
						&\le \int_0^{e^{-(\log n)^{1/4}}}(-\log t)  dt=e^{-(\log n)^{1/4}} ((\log n)^{1/4}-1),\endaligned $$
						where for the last equality we use the fact $[t(\log t-1)]'=\log t.$
These three integrals indicate 
\begin{equation}\label{II1}\int_{-\ell_1(n)}^{\ell_2(n)}x\exp\big\{-x-e^{-x}\big\}dx=\gamma+O(e^{-(\log n)^{1/4}} (\log n)^{1/4}).\end{equation}		
Using the same idea and keeping in mind that we only need to keep terms with order greater than $(\log n)^{-1},$ we get 
 \begin{equation}\label{II2} \aligned &\quad\frac{1}{4\log n}\int_{-\ell_1(n)}^{\ell_2(n)}x\exp\big\{-x-e^{-x}\big\}t_n^2(x)dx\\
 &=\frac{1}{4\log n}\int_{-\ell_1(n)}^{\ell_2(n)} x\exp\big\{-x-e^{-x}\big\}(c_n^2-2c_n x+x^2)dx\\
 &=\frac{1}{4\log n}\int_{-\infty}^{\infty}(c_n^2 x-2c_nx^2)\exp\big\{-x-e^{-x}\big\}dx +O((\log n)^{-1})\\
 &=\frac{1}{4\log n}(c_n^2\gamma -2c_n(\gamma^2+\frac{\pi^2}{6}))+O((\log n)^{-1}),\endaligned \end{equation}
 where the last equality holds since 
 $$\int_{-\infty}^{+\infty} x^2 e^{-x-e^{-x}} dx=\int_0^{\infty} (\log t)^2 e^{-t} dt=\gamma^2+\frac{\pi^2}{6} $$
 and similarly \begin{equation}\label{II3}\aligned & \quad \frac{1}{4\log n}\int_{-\ell_1(n)}^{\ell_2(n)}x\exp\big\{-2x-e^{-x}\big\} (t_n(x)^2+2t_n(x)+2)dx\\
 &=\frac{1}{4\log n} \int_{-\infty}^{\infty}(c_n^2 x-2c_n x(x+1))\exp\big\{-2x-e^{-x}\big\}dx+O((\log n)^{-1})\\
 &=\frac{1}{4\log n}(c_n^2(\gamma-1) -2c_n(\gamma^2-\gamma-1+\frac{\pi^2}{6}))+O((\log n)^{-1}).\endaligned\end{equation}				
Here, we use 
$$\aligned \int_{-\infty}^{\infty} xe^{-2x-e^{-x}}dx&=-\int_0^{\infty} t \log t e^{-t} dt=\gamma-1;  \\
		\int_{-\infty}^{\infty} x^2 e^{-2x-e^{-x}}dx&=\int_0^{\infty}  t (\log t)^2  e^{-t} dt=\gamma^2-2\gamma+\frac{\pi^2}{6}.
		 \endaligned $$	
		 The reader may find informations on these integrals via Gamma function in \cite{AS}.											 
Plugging \eqref{II1}, \eqref{II2} and \eqref{II3} back into \eqref{IIsum}, we get  
			$$\aligned
				\mathrm{II}=\gamma-\frac{c_n^2}{4\log n}+\frac{(\gamma +1)c_n}{2\log n}+O((\log n)^{-1}).
			\endaligned	$$	
		Combining the estimates for $\mathrm{I}$, $\mathrm{II}$, and $\mathrm{III},$  we finally obtain
		\begin{align*}
			\mathbb{E}[Y_n]=\gamma-\frac{c_n^2}{4\log n}+\frac{(\gamma+1) c_n}{2\log n}+O((\log n)^{-1}).
		\end{align*} 
		
	For the expectation $\mathbb{E}\left[e^{-Y_n} - \frac{t_n^2(Y_n)}{4\log n}\right]$, we follow an analogous procedure while focusing on the dominant terms. The analysis proceeds with Lemma \ref{densilem} replaced by its enhanced version \eqref{densifnew}, yielding the asymptotic expansion
	$$\aligned & \quad \mathbb{E}(e^{-Y_n}-\frac{t_n^2(Y_n)}{4\log n})\\
	&=\int_{-\ell_1(n)}^{\ell_2(n)} (e^{-x}-\frac{t_n^2(x)}{4\log n}) f_n(x) dx+O(e^{-(\log n)^{1/4}} (\log n)^{1/4})\\
	&=O((\log n)^{-9/4})+(1+o(1))\int_{-\infty}^{\infty} (e^{-x}-\frac{t_n^2(x)}{4\log n})e^{-x-e^{-x}}\bigg(1+e^{-x}\frac{t^2_n(x)+2t_n(x)+2}{4 \log n}\\
	&\quad\quad\quad\quad-\frac{t_n^2(x)}{4\log n}+\frac{t_n^4(x)}{32(\log n)^2}
	 -\frac{e^{-x}}{32(\log n)^2}\big(3t_n^4(x)+8t_n^3(x)+16t_n^2(x) +24t_n(x)+24\big)\\
	 &\quad\quad\quad\quad+\frac{e^{-2x}}{32(\log n)^2}(t_n^4(x)+4t_n^3(x)+8t_n^2(x)+8t_n(x)+4\bigg)dx.
	\endaligned $$ 
	
	To systematically analyze the integral, we decompose the integrand, denoted by \( J_n(x) \), by carefully organizing its terms according to their asymptotic orders in \( \log n \). This structured approach reveals the dominant contributions and their precise coefficients:

$$\aligned 
&\quad e^{x+e^{-x}} J_n(x)\\
&=e^{-x} +\underbrace{\sum_{k=0}^{2} \frac{A_k(x)}{\log n}}_{\text{Leading Orders}} + \underbrace{\sum_{k=0}^{4} \frac{B_k(x)}{(\log n)^2}}_{\text{Higher-Order Terms}}\\
&=e^{-x}+\frac{1}{4\log n}\big(e^{-2x}(t_n^2(x)+2t_n(x))+2)-e^{-x} t_n^2(x)-t_n^2(x)\big)\\
&\quad+\frac{1}{32(\log n)^2}\bigg(2t_n^4(x)-e^{-x}\big(t_n^4(x)+4t_n^3(x)+4t_n^2(x)\big)-e^{-2x}\big(3t_n^4(x)+8t_n^3(x)\\
&\quad+16t_n^2(x)+24t_n(x)+24\big)+e^{-3x}\big(t_n^4(x)+4t_n^3(x)+8t_n^2(x)+8t_n(x)+4)\bigg).
\endaligned $$
Now we work on the coefficient of $(\log n)^{-1}$ in $\int_{-\infty}^{+\infty}J_n(x)dx,$ that is 
 $$\aligned &\quad \int_{-\infty}^{+\infty}e^{-x-e^{-x}}\big(e^{-2x}(t_n^2(x)+2t_n(x))+2)-e^{-x} t_n^2(x)-t_n^2(x)\big)dx\\
			&=c_n^2\int_{-\infty}^{+\infty}e^{-x-e^{-x}}(e^{-2x}-e^{-x}-1))dx+2c_n\int_{-\infty}^{+\infty}e^{-x-e^{-x}}(x e^{-x}+x-e^{-2x}(x+1))\\
			&\quad+\int_{-\infty}^{+\infty}e^{-x-e^{-x}}((x^2+2x+2)e^{-2x}-x^2e^{-x}-x^2))dx\\
			&=0. \endaligned $$ 
			For the coefficients of $(\log n)^{-2}$ in $\int_{-\infty}^{+\infty}J_n(x)dx,$ we start with that of $c_n^4 (\log n)^{-2}.$ The involving integral is 
$$\aligned &\quad\int_{-\infty}^{+\infty}e^{-x-e^{-x}}(2-e^{-x}-3e^{-2x}+e^{-3x})dx&=\int_0^{+\infty} e^{-t}(2-t-3t^2+t^3)dt=1. 
\endaligned $$	 
Hence, we know 
$$ \mathbb{E}(e^{-Y_n}-\frac{t_n^2(Y_n)}{4\log n})=1+\frac{c_n^4}{32(\log n)^2}(1+o(1)).$$
The proof is then completed.
		
	\end{proof}

\begin{rmk} According to the proof of Lemma \ref{expec}, we know that 
$$\mathbb{E}|Y_n|^k=\int_{-\infty}^{\infty} e^{-e^{-x}-x} |x^k|dx+o(1)<+\infty$$
for any $k\in\mathbb{N}.$ Moreover, when necessary, we can also utilize the asymptotic of $f_n$ to derive the precise asymptotic of
$\mathbb{E}u(Y_n)$ for any measurable function $u$ on $\mathbb{R}.$
\end{rmk}
	
\section{Proof of theorem \ref{main}} 	
	With these lemmas established, we now present the proof of our main theorem. 
The following treatment of the Berry-Esseen bound was originally offered in \cite{MaMeng2025a}. 	
	
	Recall $c_n=\log\sqrt{4\pi\log n}=\frac{1}{2}\log(4\pi\log n), t_n(x)=x-c_n, $ $\ell_1(n)=\frac14\log\log n$ and $\ell_2(n)=(\log n)^{1/4}.$

\subsection{Berry-Esseen bound} 
We first prove that 
\begin{equation}\label{left}\sup_{x\le -\ell_1(n)}|F_n(x)-\exp(-e^{-x})|\ll \frac{1}{\log n}.\end{equation} 
The triangle inequality and the monotonicity  of distribution functions ensure that 
$$\sup_{x\le -\ell_1(n)}|F_n(x)-e^{-e^{-x}}|\le F_n(-\ell_1(n))+\exp(-e^{\ell_1(n)}).$$
Lemma \ref{keyFesti} indicates that 
$$F_n(-\ell_1(n))=O(\exp(-e^{\ell_1(n)}))=O(e^{-(\log n)^{1/4}}),$$
which brings the inequality \eqref{left} we want.  
Using the same tool, we get 
 \begin{equation}\label{right}\aligned \sup_{x\ge \ell_2(n)}|F_n(x)-e^{-e^{-x}}|&\le 1-F_n(\ell_2(n))+1-e^{-e^{-\ell_2(n)}}\\
 &\lesssim 1-e^{-e^{-\ell_2(n)}}\\
 &\le e^{-(\log n)^{1/4}}\ll \frac{1}{\log n}.\endaligned \end{equation}
 For the middle part, Lemma \ref{keyFesti} enables us to obtain 
  \begin{equation}\label{middle} \aligned \sup_{-\ell_1(n)\le x\le \ell_2(n)}|F_n(x)-e^{-e^{-x}}|
 &=\sup_{-\ell_1(n)\le x\le \ell_2(n)}e^{-e^{-x}} e^{-x}\frac{(1+x-c_n)^2+O(1)}{4 \log n}(1+o(1))\\
 &=\frac{c_n^2}{4e\log n}(1+o(1)),\endaligned \end{equation} 
 where the last equality holds since 
 $$\sup_{x\in\mathbb{R}}e^{-x-e^{-x}}=e^{-1}$$ and other supremums involved are bounded and the parameter $c_n\gg 1.$
Comparing these three estimates \eqref{left}, \eqref{right} and \eqref{middle}, we get 
$$\sup_{x\in\mathbb{R}}|F_n(x)-e^{-e^{-x}}|=\frac{c_n^2}{4e\log n}(1+o(1))=\frac{(\log\log n)^2}{16 e\log n}(1+o(1)).$$ 
The precise asymptotic for Berry-Esseen bound is then obtained. 

\subsection{$W_1$ Wasserstein distance} 
 
		For the $W_1$ Wasserstein distance, we use a similar approach by cutting the integral interval into three parts as  
				\begin{align*}
			W_1(F_n, \Lambda)=&\int_{|x|\ge 
		\ell_2(n)} \left|F_n(x)-e^{-e^{-x}}\right| dx+\int_{|x|< 
		\ell_2(n)} \left|F_n(x)-e^{-e^{-x}}\right| dx.
		\end{align*}
	Let $\zeta$ be a random variable obeying the Gumbel distribution. The triangle inequality and the Markov inequality work together to ensure that 	
	$$\aligned \int_{|x|\ge \ell_2(n)} |F_n(x)-e^{-e^{-x}}|&\le 
	2\int_{x\ge \ell_2(n)} (\mathbb{P}(|Y_n|\ge x)+\mathbb{P}(|\zeta|\ge x)) dx\\
	&\le 2\int_{x\ge \ell_2(n)} \frac{\mathbb{E}|Y_n|^5+\mathbb{E}|\zeta|^5}{x^5}dx\\
	&\lesssim\frac{1}{\log n}.
	\endaligned $$	
	Here, for the final inequality, we utilize the boundedness of both $\mathbb{E}|Y_n|^5$ and $\mathbb{E}|\zeta|^5$ and $\ell_2(n)=(\log n)^{1/4}.$
	Now we treat the integral $\int_{|x|< 
		\ell_2(n)} \left|F_n(x)-e^{-e^{-x}}\right| dx,$ which is continued to be cut into two parts as 
			
		\begin{equation}\label{middleinte}
			\int_{|x|< 
		\ell_2(n)} \left|F_n(x)-e^{-e^{-x}}\right| dx=(\int_{-\ell_2(n)}^{-\ell_1(n)}+\int_{-\ell_1(n)}^{\ell_2(n)})|F_n(x)-e^{-e^{-x}}|dx.
		\end{equation}
		It is clear that   
		$$\aligned \int_{-\ell_2(n)}^{-\ell_1(n)}|F_n(x)-e^{-e^{-x}}|dx&\lesssim \ell_2(n) (F_n(-\ell_1(n))+e^{-e^{\ell_1(n)}})\\
		&\lesssim (\log n)^{1/4}e^{-(\log n)^{1/4}}.
		\endaligned $$ 
We now analyze the core term. Indeed, Lemma \ref{keyFesti} yields the following estimate
		\begin{align*}
			\int_{-\ell_1(n)}^{\ell_2(n)} \left|F_n(x)-e^{-e^{-x}}\right| \, dx &=(1+o(1))\int_{-\ell_1(n)}^{\ell_2(n)}e^{-e^{-x}} e^{-x}\frac{(1+x-c_n)^2+O(1)}{4 \log n} dx.
		\end{align*}
		Employing the same approach used for the Berry-Esseen bound, we observe that the dominant term arises from the following integral 
		$$\frac{c_n^2}{4\log n}\int_{-l_1(n)}^{l_2(n)} e^{-e^{-x}} e^{-x} dx=\frac{c_n^2(1+o(1))}{4\log n}\int_{-\infty}^{+\infty} e^{-e^{-x}} e^{-x} dx=\frac{(\log\log n)^2}{16\log n}(1+o(1)).$$
		These estimates and the fact $c_n=\frac{1+o(1)}{2}\log\log n$ yield the stated result
		\begin{align*}
			W_1(F_n, \Lambda)=\frac{(\log\log n)^2}{16\log n}(1+o(1)).
		\end{align*}
	
	\subsection{Total variation} 
	By definition, 
	$$\|\mathcal{L}(Y_n)-\Lambda\|_{\rm TV}=\int_{-\infty}^{+\infty}|f_n(x)-e^{-x-e^{-x}}| dx\ge \int_{-\ell_1(n)}^{\ell_2(n)} |f_n(x)-e^{-x-e^{-x}}|dx$$ and we follow an analogous decomposition approach to get 
	$$\aligned 
	\|\mathcal{L}(Y_n)-\Lambda\|_{\rm TV}&\le \int_{-\ell_1(n)}^{\ell_2(n)} |f_n(x)-e^{-x-e^{-x}}|dx+F_n(-\ell_1(n))+e^{-e^{\ell_1(n)}}\\
	&\quad+1-F_n(\ell_2(n))+1-e^{-e^{-\ell_2(n)}}\\
	&=\int_{-\ell_1(n)}^{\ell_2(n)} |f_n(x)-e^{-x-e^{-x}}|dx+O(e^{-(\log n)^{1/4}}). 
	\endaligned $$
	For the integral above, Lemma \ref{densilem} enables us to get 
	$$\aligned &\quad \int_{-\ell_1(n)}^{\ell_2(n)} |f_n(x)-e^{-x-e^{-x}}|dx\\
	&=(1+O((\log n)^{-1}))\int_{-\ell_1(n)}^{\ell_2(n)}e^{-x-e^{-x}} |e^{-x}\frac{t_n^2(x)+2t_n(x)+2}{4 \log n}-\frac{t_n^2(x)}{4\log n}|dx.
	\endaligned $$
As before, it suffices to examine terms with \( c_n^2 \), which comes from $t_n^2(x).$
Hence,
	$$\aligned \int_{-\ell_1(n)}^{\ell_2(n)} |f_n(x)-e^{-x-e^{-x}}|dx
	&=\frac{c_n^2(1+o(1))}{4\log n}\int_{-\infty}^{+\infty}e^{-x-e^{-x}} |e^{-x}-1|dx\\
	&=\frac{c_n^2}{2e\log n}(1+o(1))=\frac{(\log\log n)^2}{8e\log n}(1+o(1)).
	\endaligned $$
	Here, we use the following integral 
	$$\aligned\int_{-\infty}^{+\infty} e^{-x-e^{-x}}|e^{-x}-1|dx&=\int_0^{+\infty} e^{-u}|u-1| du &=\frac{2}{e}.\\
	\endaligned $$
	Thus, 
	$$\|\mathcal{L}(Y_n)-\Lambda\|_{\rm TV}=\frac{(\log\log n)^2}{8e\log n}(1+o(1)).$$
		\subsection{Kullback-Leibler divergence }	
		Next, we calculate the Kullback-Leibler divergence. By definition, 
		\begin{align*} 
			D_{\rm KL}(\mathcal{L}(Y_n)|\Lambda) =\int_{-\infty}^{+\infty} f_n(x) \log \frac{f_n(x)}{e^{-x-e^{-x}}} d x. \end{align*}
			The expressions \eqref{densi} and \eqref{halffn} entail  
			\begin{equation}\label{logf}
			\log \frac{f_n(x)}{e^{-x-e^{-x}}}=e^{-x}-\frac{t_n^2(x)}{4\log n}+(n-1)\log \Phi(a_n+\frac{t_n(x)}{a_n})
		\end{equation}
		and then 
		$$D_{\rm KL}(\mathcal{L}(Y_n)|\Lambda)=\mathbb{E}\big( e^{-Y_n}-\frac{t_n^2(Y_n)}{4\log n}+(n-1)\log \Phi(a_n+\frac{Y_n-c_n}{a_n})\big). $$
		Thus, Lemma \ref{expec} helps us to get 
		\begin{align*} 
			D_{\rm KL}(\mathcal{L}(Y_n)|\Lambda) =1+\frac{c_n^4}{32(\log n)^2}(1+o(1))+(n-1)\log \Phi(a_n+\frac{Y_n-c_n}{a_n}). \end{align*}
			By definition $a_n+\frac{Y_n-c_n}{a_n}=X_{(n)},$ whose density function is $n\Phi^{n-1}(x)\phi(x)$ and then 
			\begin{align*}
			\mathbb{E}\log \Phi(a_n+\frac{Y_n-c_n}{a_n})&=\mathbb{E}\log \Phi(X_{(n)})=n\int_{-\infty}^{+\infty}	\Phi^{n-1}(x)\phi(x)\log \Phi(x) dx\\
			&=n \int_0^1 t^{n-1}\log t \; dt= -\frac{1}{n}.
			\end{align*}
Thereby, we get 
		\begin{align*}
		D_{\rm KL}(\mathcal{L}(Y_n)|\Lambda)&=\frac{c_n^4}{32(\log n)^2}(1+o(1))=\frac{(\log\log n)^4}{512(\log n)^2}(1+o(1)).
		\end{align*}
		
	\subsection{Fisher information} 
It follows from \eqref{Fishervar} that  
\begin{equation}\label{Iexpre}\aligned I(\mathcal{L}(Y_n)|\Lambda)&=\frac{1}{4}\int_{-\infty}^{+\infty}[(\log \frac{f_n(x)}{e^{-x-e^{-x}}})']^2 f_n(x)dx.
\endaligned\end{equation}
The expressions \eqref{somepoint} and  \eqref{logf} imply 
$$\aligned (\log \frac{f_n(x)}{e^{-x-e^{-x}}})'&=-e^{-x}-\frac{t_n(x)}{2\log n}+\frac{(n-1)\phi(a_n+\frac{t_n(x)}{a_n})}{a_n\Phi(a_n+\frac{t_n(x)}{a_n})}\\
&=-e^{-x}-\frac{t_n(x)}{2\log n}+\frac{(n-1)}{n\Phi(a_n+\frac{t_n(x)}{a_n})}\exp(-x-\frac{t_n^2(x)}{4\log n})\\
&=h_n(x) e^{-x}-\frac{t_n(x)}{2\log n},
\endaligned$$ 
where for simplicity we set 
$$h_n(x)=\frac{(n-1)}{n\Phi(a_n+\frac{t_n(x)}{a_n})}\exp(-\frac{t_n^2(x)}{4\log n})-1.$$
Thus, 
\begin{equation}\label{inform}\aligned  ((\log \frac{f_n(x)}{e^{-x-e^{-x}}})')^2 &=(h_n(x) e^{-x}-\frac{t_n(x)}{2\log n})^2=h^2_n(x)e^{-2x}+\frac{t_n^2(x)}{4(\log n)^2}-\frac{t_n(x)}{\log n}h_n(x)e^{-x}. 
\endaligned\end{equation}
Putting this expression back into \eqref{Iexpre}, one sees  
\begin{equation}\label{inform}\aligned I(\mathcal{L}(Y_n)|\Lambda)&=\frac14\mathbb{E} \Big(h^2_n(Y_n)e^{-2Y_n}+\frac{t_n^2(Y_n)}{4(\log n)^2}-\frac{t_n(Y_n)}{\log n}h_n(Y_n)e^{-Y_n}\Big).\endaligned\end{equation}  
As for the expectation $\mathbb{E}Y_n,$ we only need to deal with the corresponding integral over $[-\ell_1(n), \ell_2(n)].$  
Recall the density function $f_n$ of $Y_n$ satisfies  
$$f_n(x)=\exp(-x-e^{-x})\Big(1+e^{-x}\frac{t_n(x)^2+2t_n(x)+2}{4 \log n}-\frac{t_n^2(x)}{4\log n}\Big)\big(1+O((\log n)^{-1})\big)$$  
uniformly on $[-\ell_1(n), \ell_2(n)].$

The asymptotic behavior of \(I(\mathcal{L}(Y_n) \mid \Lambda)\) is determined by comparing the orders of the three terms after removing the factors \(e^{-Y_n}\) and \(e^{-2Y_n}\). This is justified because these exponential factors can be absorbed into the leading factor \(\exp(-x - e^{-x})\) of the density \(f_n\) and thus do not affect the relative growth rates of the remaining parts.  
Consequently, we only need to compare the orders of  
\begin{equation}\label{threeterm}	
h_n^2(Y_n),\qquad\frac{t_n^2(Y_n)}{(\log n)^2},\quad\text{and}\quad\frac{t_n(Y_n)h_n(Y_n)}{\log n},\end{equation}  
where the term with the highest order among these dominates the asymptotics of the whole expression.

Now on $[-\ell_1(n), \ell_2(n)],$ $h_n$ has a precise asymptotic expansion:
$$\aligned h_n(x)&=\frac{n-1}{n\Phi(a_n+\frac{t_n(x)}{a_n})}\Big(\exp\!\Big(-\frac{t_n^2(x)}{4\log n}\Big)-\frac{n\Phi(a_n+\frac{t_n(x)}{a_n})}{n-1}\Big)\\
&=\frac{n-1}{n\Phi(a_n+\frac{t_n(x)}{a_n})}\bigg(1-\frac{t_n^2(x)}{4\log n}+O\!\Big(\frac{t_n^2(x)+1}{(\log n)^2}\Big)-\Phi\!\Big(a_n+\frac{t_n(x)}{a_n}\Big)\bigg).\endaligned$$  
Using \eqref{psilem} and the relationship $\Phi(t)=1-\Psi(t),$ we obtain further  
$$\aligned
h_n(x)&=-\frac{t_n^2(x)}{4\log n}\big(1+o(1)\big)=O\!\Big(\frac{t_n^2(x)+1}{\log n}\Big)=O\big((\log n)^{-1/2}\big). \endaligned$$  
Comparing the three terms in \eqref{threeterm}, we observe that \(h_n^2(Y_n) = O\!\big((\log n)^{-1}\big)\) is of larger order than the other two terms, which are at most \(O\big((\log n)^{-5/4}\big)\). Therefore, \(h_n^2(Y_n)e^{-2Y_n}\) provides the dominant contribution to \(I(\mathcal{L}(Y_n)|\Lambda).\)

Then, using Lemma \ref{densilem}, we obtain

\[
\begin{aligned}
I(\mathcal{L}(Y_n)\mid\Lambda)
&=\frac14\bigl(1+o(1)\bigr)\int_{-\ell_1(n)}^{\ell_2(n)} h_n^2(x)\,e^{-3x-e^{-x}}
\Bigl(1+e^{-x}\frac{t_n(x)^2+2t_n(x)+2}{4\log n}-\frac{t_n^2(x)}{4\log n}\Bigr)\,dx \\[4pt]
&=\frac{c_n^4\bigl(1+o(1)\bigr)}{64(\log n)^2}\int_{-\infty}^{\infty} e^{-3x-e^{-x}}\,dx \\[4pt]
&=\frac{(\log\log n)^4}{512(\log n)^2}\,\bigl(1+o(1)\bigr).
\end{aligned}
\]
Here we have used the substitution \(t=e^{-x}\) to evaluate the integral:

\[
\int_{-\infty}^{\infty} e^{-3x-e^{-x}}\,dx
=\int_{0}^{\infty} e^{-t}\,t^{2}\,dt
=\Gamma(3)=2.
\]
Consequently,

\[
I(\mathcal{L}(Y_n)\mid\Lambda)=\frac{(\log\log n)^4}{512(\log n)^2}\,\bigl(1+o(1)\bigr).
\]
The proof is thus completed.

\section{Scaling and centering constants influence convergence rates} 
In this section, we investigate how the the scaling and centering constants influence the convergence rate. 
Recall $a_n=\sqrt{2\log n}$ and $c_n=\log(\sqrt{2\pi} a_n).$
\begin{thm}\label{main2}
Let $X_{(n)}$ be defined as above. 
\begin{itemize}
\item Let $\widetilde{b}_n$ satisfy $$2\pi \widetilde{b}_n^2 \exp(\widetilde{b}_n^2) = n^2 $$ and $\widetilde{Y}_n=\widetilde{b}_n(X_{(n)}-\widetilde{b}_n),$ we have 
\begin{equation}\label{L}\aligned &\sup_{x\in\mathbb{R}} |\mathbb{P}(X_{(n)}\le \widetilde{b}_n+\frac{x}{\widetilde{b}_n})-e^{-e^{-x}}|=\frac{d_1}{4\log n}(1+o(1));\\
&W_1(\mathcal{L}(\widetilde{Y}_n), \Lambda)=\frac{6(\gamma^2+2\gamma+2)+\pi^2}{24\log n}(1+o(1)), \endaligned \end{equation} 
where $d_1=\sup_{x\in\mathbb{R}}(x^2+2x+2) e^{-e^{-x}} e^{-x}\approx 1.3.$
\item If  
$$\widehat{b}_n=a_n-\frac{c_n}{a_n}-\frac{c_n^2-2c_n}{2a_n^3}$$ and $\widehat{Y}_n=a_n(X_{(n)}-\widehat{b}_n),$ we have 
\begin{equation}\label{M}\aligned &\sup_{x\in\mathbb{R}} |\mathbb{P}(X_{(n)}\le \widehat{b}_n+\frac{x}{a_n})-e^{-e^{-x}}|=\frac{d_2\log\log n}{4\log n}(1+o(1)); \\
& W_1(\mathcal{L}(\widehat{Y}_n), \Lambda)=\frac{(\gamma+2E_1(1))\log\log n}{4\log n}(1+o(1)) \endaligned \end{equation}	
with $d_2=\sup_{x\in\mathbb{R}}|x| e^{-e^{-x}} e^{-x}\approx 0.3$
and $\gamma+2E_1(1)=\int_0^{+\infty} e^{-t}|\log t|dt\approx 1.$\end{itemize}
\end{thm}
\begin{proof} We first prove \eqref{M}. As for the Berry-Esseen bound in Theorem \ref{main}, we only need to work on $x\in [-\ell_1(n), \ell_2(n)].$ 
By definition, 
$$
	\aligned \mathbb{P}(X_{(n)}\le \widehat{b}_n+\frac{x}{a_n})&=\mathbb{P}(X_{(n)}\le a_n+\frac{t_n(x)}{a_n}-\frac{c_n^2-2c_n}{2a_n^3})\\
&=\Phi^{n}(a_n+\frac{t_n(x)}{a_n}-\frac{c_n^2-2c_n}{2a_n^3}).\endaligned $$
 Based on the Mills ratio again  we have  
\begin{equation}\label{Psinew}\aligned & \quad \Psi(a_n+\frac{t_n(x)}{a_n}-\frac{c_n^2-2c_n}{2a_n^3})\\
&=\frac{1-(a_n+\frac{t_n(x)}{a_n}-\frac{c_n^2-2c_n}{2a_n^3})^{-2}}{\sqrt{2\pi} (a_n+\frac{t_n(x)}{a_n}-\frac{c_n^2-2c_n}{2a_n^3})}\exp(-\frac12(a_n+\frac{t_n(x)}{a_n}-\frac{c_n^2-2c_n}{2a_n^3})^2)(1+O(\frac{1}{(\log n)^2}))
\endaligned \end{equation} 
Remember $a_n\asymp \sqrt{\log n}$ and $c_n\asymp\log\log n,$ which allow us to write 
$$\aligned 1-(a_n+\frac{t_n(x)}{a_n}-\frac{c_n^2-2c_n}{2a_n^3})^{-2}&=1-(a_n+\frac{t_n(x)}{a_n})^{-2}(1+O(\frac{c_n^2}{(\log n)^2}))\\
&=(1-(a_n+\frac{t_n(x)}{a_n})^{-2})(1+O(\frac{c_n^2}{(\log n)^3})).
\endaligned $$  Via similar treatment done for other two terms in \eqref{Psinew}, we derive 
\begin{equation}\label{Psinewnew}\aligned & \quad \Psi(a_n+\frac{t_n(x)}{a_n}-\frac{c_n^2-2c_n}{2a_n^3})\\
&=\frac{(1-(a_n+\frac{t_n(x)}{a_n})^{-2})(1+O(\frac{c_n^2}{(\log n)^2}))}{\sqrt{2\pi} (a_n+\frac{t_n(x)}{a_n})}\exp(-\frac12(a_n+\frac{t_n(x)}{a_n})^2+\frac{c_n^2-2c_n}{4\log n}).
\endaligned \end{equation} 
Comparing \eqref{Psinewnew} and \eqref{Millsratio}, we have by Lemma \ref{psilem} together with 
$$\exp(\frac{c_n^2-2c_n}{4\log n})=1+\frac{c_n^2-2c_n}{4\log n}+O(\frac{c_n^4}{(\log n)^2})$$
and 
$$c_n^2-2c_n-(t_n^2(x)+2t_n(x)+2)=2c_n x-(x^2+2x+2)$$ that   
$$\aligned
\Psi(a_n+\frac{t_n(x)}{a_n}-\frac{c_n^2-2c_n}{2a_n^3})
&=\frac{1}{n}e^{-x}(1+\frac{2c_n x-(x^2+2x+2)}{4\log n})(1+O(\frac{t_n^4(x)+c_n^4}{(\log n)^{2}})).\endaligned $$
Thus, similarly as for \eqref{Fnasym1}, we derive 
$$\aligned \mathbb{P}(X_{(n)}\le \widehat{b}_n+\frac{x}{a_n})&=e^{-e^{-x}}+e^{-e^{-x}} e^{-x}\frac{-2c_n x+x^2+2x+2}{4 \log n}(1+O(\frac{(e^{-x}\vee 1)(t_n^4(x)+c_n^4)}{(\log n)^{2}}))\\
&=e^{-e^{-x}}+e^{-e^{-x}} e^{-x}\frac{-2c_n x+x^2+2x+2}{4 \log n}(1+O(\frac{1}{\log n})).
		\endaligned $$	
Consequently as for \eqref{middle}, we get 
$$\sup_{x\in\mathbb{R}}|\mathbb{P}(X_{(n)}\le \widehat{b}_n+\frac{x}{a_n})-e^{-e^{-x}}|=\frac{c_n}{2\log n}\sup_{x\in\mathbb{R}}|x| e^{-e^{-x}} e^{-x}$$ 
and then 
$$\sup_{x\in\mathbb{R}}|\mathbb{P}(X_{(n)}\le \widehat{b}_n+\frac{x}{a_n})-e^{-e^{-x}}|=\frac{d_2 \log\log n}{4\log n}(1+o(1))$$	
with $d_2=\sup_{x\in\mathbb{R}}|x| e^{-e^{-x}} e^{-x}\approx 0.2704.$

Now we prove \eqref{L}.	For $x\in [-\ell_1(n), \ell_2(n)],$ we see 
$$\aligned \Psi(\widetilde{b}_n+\frac{x}{\widetilde{b}_n})=
&=\frac{1-(\widetilde{b}_n+\frac{x}{\widetilde{b}_n})^{-2}}{\sqrt{2\pi} (\widetilde{b}_n+\frac{x}{\widetilde{b}_n})}\exp(-\frac12(\widetilde{b}_n+\frac{x}{\widetilde{b}_n})^2)(1+O(\frac{1}{(\log n)^2}))\\
&=\frac{1-\widetilde{b}_n^{-2}}{\sqrt{2\pi} \widetilde{b}_n (1+\frac{x}{\widetilde{b}_n^2})}\exp(-\frac12 \widetilde{b}_n^2-x-\frac{x^2}{2\widetilde{b}_n^2})(1+O(\frac{|x|+1}{\widetilde{b}_n^4}+\frac{1}{(\log n)^2})).
\endaligned $$
The condition on $\widetilde{b}_n$ is equivalent to $$\exp(-\frac{1}{2}\widetilde{b}_n^2)=n^{-1}\sqrt{2\pi} \;\widetilde{b}_n,$$ which ensures $\widetilde{b}_n=\sqrt{2\log n}(1+o(1)).$ This, with the fact $\frac{|x|+x^2}{\widetilde{b}_n^2}=o(1),$ yields 
\begin{equation}\label{newbn}\aligned \Psi(\widetilde{b}_n+\frac{x}{\widetilde{b}_n})&=\frac{1-\widetilde{b}_n^{-2}}{n(1+\frac{x}{\widetilde{b}_n^2})}\exp(-x-\frac{x^2}{2\widetilde{b}_n^2})(1+O(\frac{|x|+1}{(\log n)^2}))\\
&=\frac{1}{n}e^{-x}(1-\frac{x^2+2x+2}{2\widetilde{b}_n^2})(1+O(\frac{x^4+1}{(\log n)^2})) .
\endaligned \end{equation} 
Thus, similarly it follows that 
$$\aligned \mathbb{P}(X_{(n)}\le \widetilde{b}_n+\frac{x}{\widetilde{b}_n})-e^{-e^{-x}}
			&=e^{-e^{-x}} e^{-x}\frac{x^2+2x+2}{2\widetilde{b}_n^2}(1+o(1)).
		\endaligned $$ 
		Therefore,
		$$\aligned\sup_{x\in\mathbb{R}}|\mathbb{P}(X_{(n)}\le \widetilde{b}_n+\frac{x}{\widetilde{b}_n})-e^{-e^{-x}}|
		&=\frac{1+o(1)}{2\widetilde{b}_n^2}\sup_{x\in\mathbb{R}}e^{-e^{-x}} e^{-x}(x^2+2x+2)\\
		&=\frac{d_1}{4\log n}(1+o(1)).\endaligned$$ 
		Here, $$d_1=\sup_{x\in\mathbb{R}}e^{-e^{-x}} e^{-x}(x^2+2x+2)\approx 1.3.$$   
		The results on $W_1$ are obtained by integrating the corresponding functions. 
\end{proof} 

\begin{rmk}
The precise asymptotic \eqref{L} confirms the upper bound $\frac{3}{\log n}$ obtained by Hall \cite{Hall}. 	
\end{rmk}

Theorem \ref{main2} demonstrates how $b_n$ determines the convergence rate of Berry-Esseen bound as well as $W_1$ Wasserstein distance.  Next, we give convergence rates corresponding to other metrics for 
$\widehat{Y}_n=\widehat{b}_n(X_{(n)}-\widehat{b}_n).$ 
\begin{thm}\label{main3} Let $\widehat{b}_n$ and $\widehat{Y}_n$ be defined just as above. We have 
$$\aligned
			\|\mathcal{L}(\widehat{Y}_n)-\Lambda\|_{\rm TV}&=\frac{d_3(1+o(1))}{4\log n};\\
			D_{\rm KL}(\mathcal{L}(\widehat{Y}_n)|\Lambda)&=\frac{d_4}{32(\log n)^2}(1+o(1)); \\
			I(\mathcal{L}(\widehat{Y}_n)|\Lambda)&=\frac{d_5}{64(\log n)^2}(1+o(1))	
		\endaligned 	$$
		for sufficiently large $n.$  
		Here, $$\aligned d_3:&=\int_{-\infty}^{+\infty} e^{-x-e^{-x}}|-x^2+e^{-x}(x^2+2x+2)|dx\approx 2.6; \\
		d_4:&=\int_{-\infty}^{+\infty}e^{-x-e^{-x}}\big(5x^4-16x^3-4x^2+8\big)dx\approx 30.8;\\
		d_5:&=\int_{-\infty}^{+\infty}e^{-x-e^{-x}}(e^{-x}x^2+2x)^2dx\approx 11.5. \endaligned$$ 
\end{thm}
\begin{proof} The density function $\widehat{f}_n$ of $\widehat{Y}_n$ are involved in these three quantities. By definition,
$$\mathbb{P}(\widehat{Y}_n\le x)=\Phi^n(\widehat{b}_n+\frac{x}{\widehat{b}_n})$$ 
and then 
$$\widehat{f}_n(x)=\frac{n}{\sqrt{2\pi}\widehat{b}_n}\Phi^{n-1}(\widehat{b}_n+\frac{x}{\widehat{b}_n})\exp(-\frac{1}{2}(\widehat{b}_n+\frac x{\widehat{b}_n})^2)=\exp(-x-\frac{x^2}{2\widehat{b}_n^2})\Phi^{n-1}(\widehat{b}_n+\frac{x}{\widehat{b}_n}).$$ 
The same argument as for Lemma \ref{densilem}, while applying \eqref{newbn} instead of Lemma \ref{psilem}, we have 
\begin{equation}\label{newbnasy}\widehat{f}_n(x)=e^{-x-e^{-x}}(1-\frac{x^2}{2\widehat{b}_n^2}+e^{-x}\frac{x^2+2x+2}{2\widehat{b}_n^2})(1+O(\frac{x^4+1}{\widehat{b}_n^4})).\end{equation} 
Thus, 
$$\aligned \|\mathcal{L}(\widehat{Y}_n)-\Lambda\|_{\rm TV}&=\frac{1+o(1)}{2\widehat{b}_n^2}\int_{-\infty}^{+\infty} e^{-x-e^{-x}}|-x^2+e^{-x}(x^2+2x+2)|dx \\
&=\frac{d_3(1+o(1))}{4\log n}.
\endaligned $$
Now \begin{equation}\label{logfnew}\log \frac{\widehat{f}_n(x)}{e^{-x-e^{-x}}}=e^{-x}-\frac{x^2}{2\widehat{b}_n^2}+(n-1)\log \Phi(\widehat{b}_n+\frac{x}{\widehat{b}_n})\end{equation}
and then 
$$\aligned (\log \frac{\widehat{f}_n(x)}{e^{-x-e^{-x}}})'
&=-e^{-x}-\frac{x}{\widehat{b}_n^2}+\frac{(n-1)}{n\Phi(\widehat{b}_n+\frac{x}{\widehat{b}_n})}\exp(-x-\frac{x^2}{2\widehat{b}_n^2})\\
&=-\frac{1}{2\widehat{b}_n^2}(e^{-x}x^2+2x)(1+o(1))
\endaligned$$
uniformly on $[-\ell_1(n), \ell_2(n)].$
As for $I(\mathcal{L}(Y_n)|\Lambda)$, similarly we have 
$$\aligned 
I(\mathcal{L}(\widehat{Y}_n)|\Lambda)&=\frac{1+o(1)}{16\widehat{b}_n^4}\int_{-\ell_1(n)}^{\ell_2(n)}(e^{-x}x^2+2x)^2\widehat{f}_n(x) dx\\
&=\frac{1+o(1)}{16\widehat{b}_n^4}\int_{-\infty}^{+\infty}e^{-x-e^{-x}}(e^{-x}x^2+2x)^2dx\\
&=\frac{d_5(1+o(1))}{64(\log n)^2}.
\endaligned $$

As for Kullback-Leibler divergence,  we hav from \eqref{logfnew} that 
$$\aligned D_{\rm KL}(\mathcal{L}(\widehat{Y}_n)|\Lambda)&=\mathbb{E}(e^{-\widehat{Y}_n}-\frac{1}{2\widehat{b}_n^2}\mathbb{E}\widehat{Y}_n^2)-1+\frac1n.
\endaligned $$ 
We still need a new preciser version of  \eqref{newbnasy}, which is 
$$\aligned 
\widehat{f}_n(x)&=e^{-x-e^{-x}}\bigg(1+\frac{1}{2\widehat{b}_n^2}\big(e^{-x}(x^2+2x+2)-x^2)+\frac{1}{8\widehat{b}_n^4}\big(x^4+(x^2+2x+2)^2 e^{-2x}\\
&\quad\quad\quad\quad\quad-e^{-x}(x^4+4x^3+12x^2+24x+24)\big)\bigg)(1+O(\frac{(x^6+1)(e^{-x}\vee e^{-3x}\vee 1)}{(\log n)^{3}}))
\endaligned $$
for $x\in [-\ell_1(n), \ell_2(n)].$ Now $$\mathbb{E}\big(e^{-\widehat{Y}_n} - \frac{1}{2\widehat{b}_n^2}\widehat{Y}_n^2\big)=\int_{-\infty}^{+\infty} (e^{-x}-\frac{x^2}{2\widehat{b}_n^2})\widehat{f}_n(x) dx.$$ We see the integral related to  
 the remainder term $O\big(\frac{(x^6+1)(e^{-x} \vee e^{-3x} \vee 1)}{(\log n)^{3}}\big)$ in $\widehat{f}_n(x)$ is of order $(\log n)^{-3}$ because the dominated term of  $\widehat{f}_n(x)$ is $e^{-x-e^{-x}}\big(1+O(b_n^{-2}))$. This justifies the appearance of the factor \(\big(1+O((\log n)^{-3})\big)\) in the integral representation below. 
$$\aligned 
&\quad\mathbb{E}(e^{-\widehat{Y}_n}-\frac{1}{2\widehat{b}_n^2}\widehat{Y}_n^2)\\
&=(1+O((\log n)^{-3}))\int_{-\infty}^{+\infty}e^{-x-e^{-x}}(e^{-x}-\frac{x^2}{2\widehat{b}_n^2})\bigg(1+\frac{1}{2\widehat{b}_n^2}\big(e^{-x}(x^2+2x+2)-x^2)\\
&\quad+\frac{1}{8\widehat{b}_n^4}\big(x^4+(x^2+2x+2)^2 e^{-2x}-e^{-x}(x^4+4x^3+12x^2+24x+24)\big)\bigg)\\
&=1+\frac{1}{2\widehat{b}_n^2}\int_{-\infty}^{+\infty}e^{-x-e^{-x}}\left(e^{-2x}(x^2+2x+2)-x^2 e^{-x}-x^2 \right)dx\\
&\quad+\frac{1}{8\widehat{b}_n^4}\int_{-\infty}^{+\infty}e^{-x-e^{-x}}\big(2x^4-e^{-x}(x^4+4x^3+4x^2)-e^{-2x}(x^4+4x^3+12x^2+24x+24)\\
&\quad\quad\quad\quad\quad\quad+e^{-3x}(x^2+2x+2)^2\big)dx+O((\log n)^{-3}).
\endaligned $$ 
Using once the integral formula by parts, we see for any polynomial function $P(x),$ 
\begin{equation}
\int_{-\infty}^{+\infty} e^{-x-e^{-x}} e^{-(k+1)x} P(x)dx=\int_{-\infty}^{\infty} e^{-x-e^{-x}} e^{-k x} [(1+k)P(x)-P'(x)]dx	
\end{equation}
for any $k\in\mathbb{N}.$ This formula helps us to obtain 
$$\int_{-\infty}^{+\infty}e^{-x-e^{-x}}\left(e^{-2x}(x^2+2x+2)-x^2 e^{-x}-x^2 \right)dx=0$$ as well as 
$$\aligned 
&\quad\int_{-\infty}^{+\infty}e^{-x-e^{-x}}\big(2x^4-e^{-x}(x^4+4x^3+4x^2)-e^{-2x}(x^4+4x^3+12x^2+24x+24)\\
&\quad\quad +e^{-3x}(x^2+2x+2)^2\big)dx\\
&=\int_{-\infty}^{+\infty}e^{-x-e^{-x}}\big(5x^4-16x^3-4x^2+8\big)dx\\
&\approx 30.8.
\endaligned $$
This yields 
$$D_{\rm KL}(\mathcal{L}(\widehat{Y}_n)|\Lambda)=\frac{d_4(1+o(1))}{8\widehat{b}_n^4}=\frac{d_4(1+o(1))}{32(\log n)^2}.$$ 
The proof is then completed.  

\begin{rmk}
Let $$\widetilde{b}_n=a_n-\frac{c_n}{a_n}-\frac{c_n^2-2c_n}{2a_n^3}$$ and define $\widetilde{Y}_n=a_n(X_{(n)}-\widetilde{b}_n).$ From the proof of Theorem \ref{main3}, we deduce the following asymptotic results for sufficiently large $n$  
$$\|\mathcal{L}(\widetilde{Y}_n)-\Lambda\|_{\rm TV}=\frac{d_6\log\log n}{\log n}(1+o(1))$$ and 
$$ 
D_{\rm KL}(\mathcal{L}(\widetilde{Y}_n)|\Lambda)=\frac{d_7(\log\log n)^2}{(\log n)^2}(1+o(1))\quad  \text{and} \quad 
			I(\mathcal{L}(\widetilde{Y}_n)|\Lambda)=\frac{d_8(\log\log n)^2}{(\log n)^2}(1+o(1)),
$$  
where $d_6, d_7, d_8$ are the suprema of certain explicitly given functions. 
\end{rmk}
\end{proof} 

%
%
%
%
%
%
%
%
%
%

\end{document}